\documentclass[12pt,vatola]{article}

\usepackage{graphicx}

\def\c{\centerline}

\def\re#1{\par\hangindent\parindent\indent\llap{#1\enspace}\ignorespaces}

\def\no{\noindent}

\begin{document}

\vskip 10mm

\c{\bf\large Parallel Bundles in Planar Map Geometries}

\vskip 8mm

\c{Linfan Mao}\vskip 5mm

\c{\scriptsize Institute of Systems Science of Academy of
Mathematics and Systems}

\c{\scriptsize Chinese Academy of Sciences, Beijing 100080,
P.R.China}

\c{\scriptsize E-mail: maolinfan@163.com}

\vskip 5mm

\begin{minipage}{130mm}

\no{\bf Abstract}: {\small Parallel lines are very important
objects in Euclid plane geometry and its behaviors can be gotten
by one's intuition. But in a planar map geometry, a kind of the
Smarandache geometries, the situation is complex since it may
contains elliptic or hyperbolic points. This paper concentrates on
the behaviors of parallel bundles in planar map geometries, a
generalization of parallel lines in plane geometry and obtains
characteristics for parallel bundles.}

\no{\bf Key Words}: {\small parallel bundle, planar map,
Smarandache geometry, map geometry, classification.}

\no{\bf AMS(2000)}: {\small 05C15, 20H15, 51D99, 51M05}

\end{minipage}

\vskip 6mm

{\bf $1$ Introduction}

\vskip 3mm

A {\it map} is a connected topological graph cellularly embedded
in a surface. On the past century, many works are concentrated on
to find the combinatorial properties of maps, such as to determine
whether exists a particularly embedding on a surface ($[7][11]$)
or to enumerate a family of maps ($[6]$). All these works are on
the side of algebra, not the object itself, i.e., geometry. For
the later, more attentions are given to its element's behaviors,
such as, the line, angle, area, curvature, $\cdots$, see also
$[12]$ and $[14]$. For returning to its original face, the
conception of map geometries is introduced in $[10]$. It is proved
in $[10]$ that the map geometries are nice model of the
Smarandache geometries. They are also a new kind of intrinsic
geometry of surfaces ($[1]$). The main purpose of this paper is to
determine the behaviors of parallel bundles in planar geometries,
a generalization of parallel lines in the Euclid plane geometry.

An axiom is said {\it Smarandachely denied} if the axiom behaves
in at least two different ways within the same space, i.e.,
validated and invalided, or only invalided but in multiple
distinct ways.

A {\it Smarandache geometry} is a geometry which has at least one
Smarandachely denied axiom($1969$)($[5][13]$).

In $[3][4]$, Iseri presented a nice model of the Smarandache
geometries, called $s$-manifolds by using equilateral triangles,
which is defined as follows($[3][5][9]$):

{\it An $s$-manifold is any collection ${\mathcal C}(T,n)$ of
these equilateral triangular disks $T_i, 1\leq i\leq n$ satisfying
the following conditions:}

$(i)$ {\it Each edge $e$ is the identification of at most two
edges $e_i,e_j$ in two distinct triangular disks $T_i,T_j, 1\leq
i,j\leq n$ and $i\not= j$;}

$(ii)$ {\it Each vertex $v$ is the identification of one vertex in
each of five, six or seven distinct triangular disks.}

The conception of map geometries without boundary is defined as
follows ($[10]$).

\vskip 3mm

\no{\bf Definition $1.1$} {\it For a given combinatorial map $M$,
associates a real number $\mu (u), 0 \ < \mu (u) \ < \ \pi$, to
each vertex $u, u\in V(M)$. Call $(M,\mu)$ a map geometry without
boundary, $\mu (u)$ the angle factor of the vertex $u$ and to be
orientablle or non-orientable if $M$ is orientable or not.}

\vskip 2mm

In $[10]$, it has proved that map geometries are the Smarandache
geometries. The realization of each vertex $u, u\in V(M)$ in $R^3$
space is shown in the Fig.$1$ for each case of $\rho (u)\mu (u) \
> \ 2\pi$, $=2\pi$ or $ \ < \ 2\pi$, call {\it elliptic point,
euclidean point} and {\it hyperbolic point}, respectively.

\includegraphics[bb=20 20 100 140]{g4.eps}

\vskip 2mm

\c{$\rho (u)\mu (u) \ < \ 2\pi$\hskip 15mm $\rho (u)\mu (u)
=2\pi$\hskip 20mm $\rho (u)\mu (u) \ > \ 2\pi$}\vskip 2mm

\c{\bf Fig.$1$}\vskip 2mm

\no Therefore, a line passes through an elliptic vertex, an
euclidean vertex or a hyperbolic vertex $u$ has angle $\frac{\rho
(u)\mu (u)}{2}$ at the vertex $u$. It is not $180^{\circ}$ if the
vertex $u$ is elliptic or hyperbolic. Then what is the angle of a
line passes through a point on an edge of a map? It is
$180^{\circ}$? Since we wish the change of angles on an edge is
smooth, the answer is not. For the Smarandache geometries, the
parallel lines in them are need to be given more attention. We
have the following definition.

\vskip 3mm

\no{\bf Definition $1.2$} {\it A family $\mathcal L$ of infinite
lines not intersecting each other in a planar geometry is called a
parallel bundle.}

\vskip 2mm

In the Fig.$2$, we present all cases of parallel bundles passing
through an edge in planar geometries, where, (a) is the case of
points $u,v$ are same type with $\rho (u)\mu (u)=\rho (v)\mu (v)$,
(b) and (c) the cases of same types with $\rho (u)\mu (u) \ > \
\rho (v)\mu (v)$ and (d) the case of $u$ is elliptic and $v$
hyperbolic.

\vskip 3mm

\includegraphics[bb=5 5 10 140]{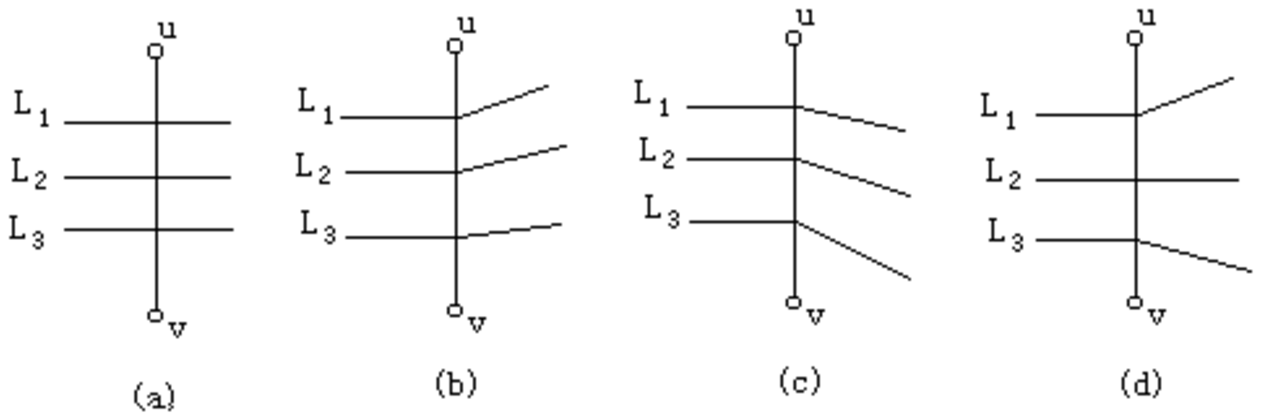}

\vskip 2mm

\c{\bf Fig.$2$}

\no Here, we assume the angle at the intersection point is in
clockwise, that is, a line passing through an elliptic point will
bend up and a hyperbolic point will bend down, such as the cases
(b),(c) in the Fig.$2$. For a vector $\overrightarrow{O}$ on the
Euclid plane, call it an {\it orientation}. We classify parallel
bundles in planar map geometries along an orientation
$\overrightarrow{O}$.

\vskip 5mm

{\bf $2.$ A condition for parallel bundles}

\vskip 3mm

We investigate the behaviors of parallel bundles in the planar map
geometries. For this object, we define a function $f(x)$ of angles
on an edge of a planar map as follows.

\vskip 3mm

\no{\bf Definition $2.1$} {\it Denote by $f(x)$ the angle function
of a line $L$ passing through an edge $uv$ at the point of
distance $x$ to $u$ on the edge $uv$.}

\vskip 2mm

Then we get the following result.

\vskip 3mm

\no{\bf Proposition $2.1$} {\it A family $\mathcal L$ of parallel
lines passing through an edge $uv$ is a parallel bundle iff }

$$\left. \frac{df}{dx} \right|_+\geq 0.$$

\vskip 2mm

{\it Proof} If $\mathcal L$ is a parallel bundle, then any two
lines $L_1,L_2$ will not intersect after them passing through the
edge $uv$. Therefore, if $\theta_1, \theta_2$ are the angles of
$L_1, L_2$ at the intersect points of $L_1, L_2$ with $uv$ and
$L_2$ is far from $u$ than $L_2$, then we know that $\theta_2\geq
\theta_1$. Whence, for any point with $x$ distance from $u$ and
$\Delta x \ > \ 0$, we have that

$$f(x+\Delta x)-f(x)\geq 0.$$

\no Therefore, we get that

$$\left. \frac{df}{dx} \right|_+=\lim_{\Delta x\to +0}\frac{f(x+\Delta x)-f(x)}
{\Delta x}\geq 0.$$

\no As the cases in the Fig.$1$.

Now if $\left.\frac{df}{dx}\right|_+\geq 0$, then $f(y)\geq f(x)$
if $y\geq x$. Since ${\mathcal L}$ is a family of parallel lines
before meeting $uv$, whence, any two lines in ${\mathcal L}$ will
not intersect each other after them passing through $uv$.
Therefore, ${\mathcal L}$ is a parallel bundle.\quad\quad
$\natural$

A general condition for a family of parallel lines passing through
a cut of a planar map being a parallel bundle is the following.

\vskip 3mm

\no{\bf Proposition $2.2$} {\it Let $(M,\mu)$ be a planar map
geometry, $C=\{u_1v_1,u_2v_2,\cdots ,u_lv_l\}$ a cut of the map
$M$ with order $u_1v_1,u_2v_2,\cdots ,u_lv_l$ from the left to the
right, $l\geq 1$ and the angle functions on them are $f_1,
f_2,\cdots ,f_l$, respectively, also see the Fig.$3$.

\includegraphics[bb=5 5 10 140]{h2.eps}

\vskip 2mm

\c{\bf Fig.$3$}

\no Then a family $\mathcal L$ of parallel lines passing through
$C$ is a parallel bundle iff for any $x, x \geq 0$,}

\begin{eqnarray*}
f'_1(x)\geq 0\\
f'_{1+}(x)+f'_{2+}(x)\geq 0\\
f'_{1+}(x)+f'_{2+}(x)+f'_{3+}(x)\geq 0\\
\cdots\cdots\cdots\cdots\\
f'_{1+}(x)+f'_{2+}(x)+\cdots +f'_{l+}(x)\geq 0.
\end{eqnarray*}

\vskip 2mm

{\it Proof} According to the Proposition $2.1$, see the following
Fig.$4$,

\includegraphics[bb=5 5 10 140]{h3.eps}

\vskip 2mm

\c{\bf Fig.$4$}

\no we know that any lines will not intersect after them passing
through $u_1v_1$ and $u_2v_2$  iff for $\forall \Delta x \ > 0$
and $x \geq 0$,

$$f_2(x+\Delta x)+f'_{1+}(x)\Delta x\geq f_2(x).$$

\no That is,

$$f'_{1+}(x)+f'_{2+}(x)\geq 0.$$

Similarly, any lines will not intersect after them passing through
$u_1v_1, u_2v_2$ and $u_3v_3$ iff for $\forall \Delta x \ > 0$ and
$x \geq 0$,

$$f_3(x+\Delta x)+f'_{2+}(x)\Delta x+f'_{1+}(x)\Delta x\geq f_3(x).$$

\no That is,

$$f'_{1+}(x)+f'_{2+}(x)+f'_{3+}(x)\geq 0.$$

\no Generally, any lines will not intersect after them passing
through $u_1v_1, u_2v_2, \cdots , u_{l-1}v_{l-1}$ and $u_lv_l$ iff
for $\forall \Delta x \ > 0$ and $x \geq 0$,

$$f_l(x+\Delta x)+f'_{l-1+}(x)\Delta x+\cdots +f'_{1+}(x)\Delta x\geq f_l(x).$$

\no Whence, we get that

$$f'_{1+}(x)+f'_{2+}(x)+\cdots +f'_{l+}(x)\geq 0.$$

Therefore, a family $\mathcal L$ of parallel lines passing through
$C$ is a parallel bundle iff for any $x, x \geq 0$, we have that

\begin{eqnarray*}
f'_1(x)\geq 0\\
f'_{1+}(x)+f'_{2+}(x)\geq 0\\
f'_{1+}(x)+f'_{2+}(x)+f'_{3+}(x)\geq 0\\
\cdots\cdots\cdots\cdots\\
f'_{1+}(x)+f'_{2+}(x)+\cdots +f'_{l+}(x)\geq 0.
\end{eqnarray*}

\no This completes the proof. \quad\quad $\natural$.

\vskip 3mm

\no{\bf Corollary $2.1$} {\it Let $(M,\mu)$ be a planar map
geometry, $C=\{u_1v_1,u_2v_2,\cdots ,u_lv_l\}$ a cut of the map
$M$ with order $u_1v_1,u_2v_2,\cdots ,u_lv_l$ from the left to the
right, $l\geq 1$ and the angle functions on them are $f_1,
f_2,\cdots ,f_l$. Then a family $\mathcal L$ of parallel lines
passing through $C$ is still parallel lines after them leaving $C$
iff for any $x, x \geq 0$,}

\begin{eqnarray*}
f'_1(x)\geq 0\\
f'_{1+}(x)+f'_{2+}(x)\geq 0\\
f'_{1+}(x)+f'_{2+}(x)+f'_{3+}(x)\geq 0\\
\cdots\cdots\cdots\cdots\\
f'_{1+}(x)+f'_{2+}(x)+\cdots +f'_{l-1+}(x)\geq 0\\
f'_{1+}(x)+f'_{2+}(x)+\cdots +f'_{l+}(x)=0.
\end{eqnarray*}

\vskip 2mm

{\it Proof} According to the Proposition $2.2$, we know the
condition is a necessary and sufficient condition for $\mathcal L$
being a parallel bundle. Now since lines in $\mathcal L$ are
parallel lines after them leaving $C$ iff for any $x\geq 0$ and
$\Delta x\geq 0$, there must be that

$$f_l(x+\Delta x)+f'_{l-1+}(x)\Delta x+\cdots +f'_{1+}(x)\Delta x= f_l(x).$$

\no Therefore, we get that

$$f'_{1+}(x)+f'_{2+}(x)+\cdots +f'_{l+}(x)=0\quad\quad \natural$$

When do the parallel lines parallel the initial parallel lines
after them passing through a cut $C$ in a planar map geometry? The
answer is in the following result.

\vskip 3mm

\no{\bf Proposition $2.3$} {\it Let $(M,\mu)$ be a planar map
geometry, $C=\{u_1v_1,u_2v_2,\cdots ,u_lv_l\}$ a cut of the map
$M$ with order $u_1v_1,u_2v_2,\cdots ,u_lv_l$ from the left to the
right, $l\geq 1$ and the angle functions on them are $f_1,
f_2,\cdots ,f_l$. Then the parallel lines parallel the initial
parallel lines after them passing through $C$ iff for $\forall
x\geq 0$,}

\begin{eqnarray*}
f'_1(x)\geq 0\\
f'_{1+}(x)+f'_{2+}(x)\geq 0\\
f'_{1+}(x)+f'_{2+}(x)+f'_{3+}(x)\geq 0\\
\cdots\cdots\cdots\cdots\\
f'_{1+}(x)+f'_{2+}(x)+\cdots +f'_{l-1+}(x)\geq 0
\end{eqnarray*}

\no{\it and}

$$f_1(x)+f_2(x)+\cdots +f_l(x)=l\pi.$$

\vskip 2mm

{\it Proof} According to the Proposition $2.2$ and Corollary
$2.1$, we know the parallel lines passing through $C$ is a
parallel bundle.

We calculate the angle $\alpha (i,x)$ of a line $L$ passing
through an edge $u_iv_i, 1\leq i\leq l$ with the line before it
meeting $C$ at the intersection of $L$ with the edge $u_iv_i$,
where $x$ is the distance of the intersection point to $u_1$ on
$u_1v_1$, see also the Fig.$4$. By the definition, we know the
angle $\alpha (1,x)=f(x)$ and $\alpha
(2,x)=f_2(x)-(\pi-f_1(x))=f_1(x)+f_2(x)-\pi$.

Now if $\alpha (i,x)=f_1(x)+f_2(x)+\cdots +f_i(x)-(i-1)\pi$, then
similar to the case $i=2$, we know that $\alpha
(i+1,x)=f_{i+1}(x)-(\pi -\alpha (i,x))=f_{i+1}(x)+\alpha
(i,x)-\pi$. Whence, we get that

$$\alpha (i+1,x)=f_1(x)+f_2(x)+\cdots +f_{i+1}(x)-i\pi.$$

Notice that a line $L$ parallel the initial parallel line after it
passing through $C$ iff $\alpha (l,x)=\pi$, i.e.,

$$f_1(x)+f_2(x)+\cdots +f_l(x)=l\pi.$$

This completes the proof. \quad\quad $\natural$

\vskip 5mm

{\bf $3.$ Linear condition and combinatorial realization for
parallel bundles}

\vskip 3mm

For the simplicity, we can assume the function $f(x)$ is linear
and denoted it by $f_l(x)$. We can calculate $f_l(x)$ as follows.

\vskip 3mm

\no{\bf Proposition $3.1$} {\it The angle function $f_l(x)$ of  a
line $L$ passing through an edge $uv$ at the point with distance
$x$ to u is }

 $$f_l(x)=
(1-\frac{x}{d(uv)})\frac{\rho (u)\mu
(v)}{2}+\frac{x}{d(uv)}\frac{\rho (v)\mu (v)}{2},$$

\no{\it where, $d(uv)$ is the length of the edge $uv$.}

\vskip 2mm

{\it Proof} Since $f_l(x)$ is linear, we know that $f_l(x)$
satisfies the following equation.

$$\frac{f_l(x)-\frac{\rho (u)\mu (u)}{2}}{\frac{\rho (v)\mu (v)}{2}-\frac{\rho (u)\mu (u)}{2}}=\frac{x}{d(uv)},$$

\no Calculation shows that

 $$f_l(x)=
(1-\frac{x}{d(uv)})\frac{\rho (u)\mu
(v)}{2}+\frac{x}{d(uv)}\frac{\rho (v)\mu (v)}{2}.\quad \natural$$

\vskip 3mm

\no{\bf Corollary $3.1$} {\it Under the linear assumption, a
family $\mathcal L$ of parallel lines passing through an edge $uv$
is a parallel bundle iff}

$$\frac{\rho (u)}{\rho (v)}\leq\frac{\mu (v)}{\mu (u)}.$$

\vskip 2mm

{\it Proof} According to the Proposition $2.1$, a family of
parallel lines passing through an edge $uv$ is a parallel bundle
iff for $\forall x, x \geq 0$, $f'(x)\geq 0$, i.e.,

$$\frac{\rho (v)\mu (v)}{2d(uv)}-\frac{\rho (u)\mu (u)}{2d(uv)}\geq 0.$$

Therefore, a family $\mathcal L$ of parallel lines passing through
an edge $uv$ is a parallel bundle iff

$$ \rho (v)\mu (v)\geq \rho (u)\mu (u).$$

\no Whence,

$$\frac{\rho (u)}{\rho (v)}\leq\frac{\mu (v)}{\mu (u)}.\quad\quad \natural$$

For a family of parallel lines pass through a cut, we have the
following condition for it being a parallel bundle.

\vskip 3mm

\no{\bf Proposition $3.2$} {\it Let $(M,\mu)$ be a planar map
geometry, $C=\{u_1v_1,u_2v_2,\cdots ,u_lv_l\}$ a cut of the map
$M$ with order $u_1v_1,u_2v_2,\cdots ,u_lv_l$ from the left to the
right, $l\geq 1$. Then under the linear assumption, a family $L$
of parallel lines passing through $C$ is a parallel bundle iff the
angle factor $\mu$ satisfies the following linear inequality
system}

$$\rho (v_1)\mu (v_1)\geq \rho (u_1)\mu (u_1)$$

$$\frac{\rho (v_1)\mu (v_1)}{d(u_1v_1)}+\frac{\rho (v_2)\mu
(v_2)}{d(u_2v_2)}\geq \frac{\rho (u_1)\mu
(u_1)}{d(u_1v_1)}+\frac{\rho (u_2)\mu (u_2)}{d(u_2v_2)}$$

$$\cdots\cdots\cdots\cdots$$

\begin{eqnarray*}
\frac{\rho (v_1)\mu (v_1)}{d(u_1v_1)}&+&\frac{\rho (v_2)\mu
(v_2)}{d(u_2v_2)}+\cdots +\frac{\rho(v_l)\mu (v_l)}{d(u_lv_l)}\\
&\geq &\frac{\rho (u_1)\mu (u_1)}{d(u_1,v_1)}+\frac{\rho (u_2)\mu
(u_2)}{d(u_2,v_2)}+\cdots + \frac{\rho (u_l)\mu
(u_l)}{d(u_l,v_l)}.\end{eqnarray*}

{\it Proof} Under the linear assumption, for any integer $i, i\geq
1$, we know that

$$f'_{i+}(x)=\frac{\rho (v_i)\mu (v_i)-\rho (u_i)\mu (u_i)}{2d(u_iv_i)}$$

\no by the Proposition $3.1$. Whence, according to the Proposition
$2.2$, we get that a family $L$ of parallel lines passing through
$C$ is a parallel bundle iff the angle factor $\mu$ satisfies the
following linear inequality system

$$\rho (v_1)\mu (v_1)\geq \rho (u_1)\mu (u_1)$$

$$\frac{\rho (v_1)\mu (v_1)}{d(u_1v_1)}+\frac{\rho (v_2)\mu
(v_2)}{d(u_2v_2)}\geq \frac{\rho (u_1)\mu
(u_1)}{d(u_1v_1)}+\frac{\rho (u_2)\mu (u_2)}{d(u_2v_2)}$$

$$\cdots\cdots\cdots\cdots$$

\begin{eqnarray*}
\frac{\rho (v_1)\mu (v_1)}{d(u_1v_1)}&+&\frac{\rho (v_2)\mu
(v_2)}{d(u_2v_2)}+\cdots +\frac{\rho(v_l)\mu (v_l)}{d(u_lv_l)}\\
&\geq &\frac{\rho (u_1)\mu (u_1)}{d(u_1,v_1)}+\frac{\rho (u_2)\mu
(u_2)}{d(u_2,v_2)}+\cdots + \frac{\rho (u_l)\mu
(u_l)}{d(u_l,v_l)}.\end{eqnarray*}

\no This completes the proof. \quad\quad $\natural$

For planar maps underlying a regular graph, we have the following
interesting results for parallel bundles.

\vskip 3mm

\no{\bf Corollary $3.2$}  {\it Let $(M,\mu)$ be a planar map
geometry with $M$ underlying a regular graph,
$C=\{u_1v_1,u_2v_2,\cdots ,u_lv_l\}$ a cut of the map $M$ with
order $u_1v_1,u_2v_2,\cdots ,u_lv_l$ from the left to the right,
$l\geq 1$. Then under the linear assumption, a family $L$ of
parallel lines passing through $C$ is a parallel bundle iff the
angle factor $\mu$ satisfies the following linear inequality
system}

$$\mu (v_1)\geq  \mu (u_1)$$

$$\frac{\mu (v_1)}{d(u_1v_1)}+\frac{\mu (v_2)}{d(u_2v_2)}\geq
\frac{\mu (u_1)}{d(u_1v_1)}+\frac{\mu (u_2)}{d(u_2v_2)}$$

$$\cdots\cdots\cdots\cdots$$

$$ \frac{\mu (v_1)}{d(u_1v_1)}+\frac{\mu (v_2)}{d(u_2v_2)}+\cdots
+\frac{\mu (v_l)}{d(u_lv_l)}\geq \frac{\mu
(u_1)}{d(u_1v_1)}+\frac{\mu (u_2)}{d(u_2v_2)}+\cdots + \frac{\mu
(u_l)}{d(u_lv_l)}$$

\no{\it and particularly, if assume that all the lengths of edges
in $C$ are the same, then }

\begin{eqnarray*}
\mu (v_1)&\geq& \mu (u_1)\\
\mu (v_1)+\mu (v_2)&\geq& \mu (u_1)+\mu (u_2)\\
\cdots\cdots &\cdots& \cdots\cdots\\
\mu (v_1)+\mu (v_2)+\cdots +\mu (v_l)&\geq& \mu (u_1)+\mu
(u_2)+\cdots +\mu (u_l).
\end{eqnarray*}

Certainly, by choosing different angle factors, we can also get
combinatorial conditions for existing parallel bundles under the
linear assumption.

\vskip 3mm

\no{\bf Proposition $3.3$} {\it Let $(M,\mu)$ be a planar map
geometry, $C=\{u_1v_1,u_2v_2,\cdots ,u_lv_l\}$ a cut of the map
$M$ with order $u_1v_1,u_2v_2,\cdots ,u_lv_l$ from the left to the
right, $l\geq 1$. If for any integer $i,i\geq 1$,}

$$\frac{\rho (u_i)}{\rho (v_i)}\leq\frac{\mu (v_i)}{\mu (u_i)},$$

\no{then under the linear assumption, a family $L$ of parallel
lines passing through $C$ is a parallel bundle.}

\vskip 2mm

{\it Proof} Notice that under the linear assumption, for any
integer $i, i\geq 1$, we know that

$$f'_{i+}(x)=\frac{\rho (v_i)\mu (v_i)-\rho (u_i)\mu (u_i)}{2d(u_iv_i)}$$

\no by the Proposition $3.1$. Whence, $f'_{i+}(x)\geq 0$ for
$i=1,2,\cdots ,l$. Therefore, we get that

\begin{eqnarray*}
f'_1(x)\geq 0\\
f'_{1+}(x)+f'_{2+}(x)\geq 0\\
f'_{1+}(x)+f'_{2+}(x)+f'_{3+}(x)\geq 0\\
\cdots\cdots\cdots\cdots\\
f'_{1+}(x)+f'_{2+}(x)+\cdots +f'_{l+}(x)\geq 0.
\end{eqnarray*}

By the Proposition $2.2$, we know that a family $L$ of parallel
lines passing through $C$ is a parallel bundle.\quad\quad
$\natural$

\vskip 5mm

{\bf $4.$ Classification of parallel bundles}

\vskip 3mm

For a cut $C$ in a planar map geometry and $e\in C$, denote by
$f_e(x)$ the angle function on the edge $e$,
$f(C,x)=\sum\limits_{e\in C}f_e(x)$. If $f(C,x)$ is independent on
$x$, then we abbreviate it to $f(C)$. According to the results in
the Section $2$ and $3$, we can classify the parallel bundles with
a given orientation $\overrightarrow{O}$ in planar map geometries
into the following $15$ classes, where, each class is labelled by
a $4$-tuple $0,1$ code.

\vskip 3mm

{\bf  Classification of parallel bundles}\vskip 2mm

$(1)$ \ {\bf $\mathcal{C}_{1000}$}: {\it for any cut $C$ along
$\overrightarrow{O}$, $f(C)=|C|\pi$;}

$(2)$ \ {\bf $\mathcal{C}_{0100}$}: {\it for any cut $C$ along
$\overrightarrow{O}$, $f(C) \ < \ |C|\pi$;}

$(3)$ \ {\bf $\mathcal{C}_{0010}$}: {\it for any cut $C$ along
$\overrightarrow{O}$, $f(C) \ > \ |C|\pi$ ;}

$(4)$ \ {\bf $\mathcal{C}_{0001}$}: {\it for any cut $C$ along
$\overrightarrow{O}$, $f'_+(C,x) \ > \ 0$ for $\forall x, x \geq
0;$}

$(5)$ \ {\bf $\mathcal{C}_{1100}$}: {\it There exist cuts
$C_1,C_2$ along $\overrightarrow{O}$, such that $f(C_1)=|C_1|\pi$
and $f(C_2)=c \ < \ |C_2|\pi$;}

$(6)$ \ {\bf $\mathcal{C}_{1010}$}: {\it there exist cuts
$C_1,C_2$ along $\overrightarrow{O}$, such that $f(C_1)=|C_1|\pi$
and $f(C_2) \ > \ |C_2|\pi$;}

$(7)$ \ {\bf $\mathcal{C}_{1001}$}: {\it there exist cuts
$C_1,C_2$ along $\overrightarrow{O}$, such that $f(C_1)=|C_1|\pi$
and $f'_+(C_2,x) \ > \ 0$ for $\forall x, x \geq 0;$}

$(8)$ \ {\bf $\mathcal{C}_{0110}$}: {\it there exist cuts
$C_1,C_2$ along $\overrightarrow{O}$, such that $f(C_1) \ < \
|C_1|\pi$ and $f(C_2) \ > \ |C_2|\pi$;}

$(9)$ \ {\bf $\mathcal{C}_{0101}$}: {\it there exist cuts
$C_1,C_2$ along $\overrightarrow{O}$, such that $f(C_1) \ < \
|C_1|\pi$ and $f'_+(C_2,x) \ > \ 0$ for $\forall x, x \geq 0;$}

$(10)$ \ {\bf $\mathcal{C}_{0011}$}: {\it there exist cuts
$C_1,C_2$ along $\overrightarrow{O}$, such that $f(C_1) \ > \
|C_1|\pi$ and $f'_+(C_2,x) \ > \ 0$ for $\forall x, x \geq 0;$}

$(11)$ \ {\bf $\mathcal{C}_{1110}$}: {\it there exist cuts
$C_1,C_2$ and $C_3$ along $\overrightarrow{O}$, such that $f(C_1)=
|C_1|\pi$, $f(C_2) \ < \ |C_2|\pi$ and $f(C_3) \ > \ |C_3|\pi$;}

$(12)$ \ {\bf $\mathcal{C}_{1101}$}: {\it there exist cuts
$C_1,C_2$ and $C_3$ along $\overrightarrow{O}$, such that $f(C_1)=
|C_1|\pi$, $f(C_2) \ < \ |C_2|\pi$ and $f'_+(C_3,x) \
> \ 0$ for $\forall x, x \geq 0;$}

$(13)$ \ {\bf $\mathcal{C}_{1011}$}: {\it there exist cuts
$C_1,C_2$ and $C_3$ along $\overrightarrow{O}$, such that $f(C_1)=
|C_1|\pi$, $f(C_2) \ > \ |C_2|\pi$ and $f'_+(C_1,x) \
> \ 0$ for $\forall x, x \geq 0;$}

$(14)$ \ {\bf $\mathcal{C}_{0111}$}: {\it there exist cuts
$C_1,C_2$ and $C_3$ along $\overrightarrow{O}$, such that $f(C_1)
\ < \ |C_1|\pi$, $f(C_2) \ > \ |C_2|\pi$ and $f'_+(C_1,x) \ > \ 0$
for $\forall x, x \geq 0;$}

$(15)$ \ {\bf $\mathcal{C}_{1111}$}: {\it there exist cuts
$C_1,C_2, C_3$ and $C_4$ along $\overrightarrow{O}$, such that
$f(C_1)= |C_1|\pi$, $f(C_2) \ < \ |C_2|\pi$, $f(C_3) \ > \
|C_3|\pi$ and $f'_+(C_4,x) \
> \ 0$ for $\forall x, x \geq 0.$}

\vskip 2mm

Notice that only the first three classes may be parallel lines
after them passing through the cut $C$. All of the other classes
are only parallel bundles, not parallel lines in the usual
meaning.

\vskip 3mm

\no{\bf Proposition $4.1$} {\it For an orientation
$\overrightarrow{O}$, the $15$ classes
$\mathcal{C}_{1000}\sim\mathcal{C}_{1111}$ are all the parallel
bundles in planar map geometries.}

\vskip 2mm

{\it Proof} Not loss of generality, we assume $C_1,C_2,\cdots
,C_m, m\geq 1$, are all the cuts along $\overrightarrow{O}$ in a
planar map geometry $(M,\mu )$ from the upon side of
$\overrightarrow{O}$ to its down side. We find their structural
characters for each case in the following discussion.

$\mathcal{C}_{1000}$: \ By the Proposition $2.3$, a family
${\mathcal L}$ of parallel lines parallel their initial lines
before meeting $M$ after the passing through $M$.

$\mathcal{C}_{0100}$: \ By the definition, a family ${\mathcal L}$
of parallel lines is a parallel bundle along $\overrightarrow{O}$
only if

$$f(C_1)\leq f(C_2)\leq \cdots \leq f(C_m) \ < \pi.$$

\no Otherwise, some lines in ${\mathcal L}$ will intersect.
According to the Corollary $2.1$, they parallel each other after
they passing through $M$ only if

$$f(C_1)= f(C_2)= \cdots = f(C_m) \ < \pi.$$

$\mathcal{C}_{0010}$: Similar to the case $\mathcal{C}_{0100}$, a
family ${\mathcal L}$ of parallel lines is a parallel bundle along
$\overrightarrow{O}$ only if

$$\pi \ < \ f(C_1)\leq f(C_2)\leq \cdots \leq f(C_m)$$

\no and parallel each other after they passing through $M$ only if

$$\pi \ < \ f(C_1)= f(C_2)= \cdots = f(C_m).$$

$\mathcal{C}_{0001}$: Notice that by the proof of the Proposition
$2.3$, a line has angle $f(C,x)-(|C|-1)\pi$ after it passing
through $C$ with the initial line before meeting $C$. In this
case, a family ${\mathcal L}$ of parallel lines is a parallel
bundle along $\overrightarrow{O}$ only if for $\forall x_i,
x_i\geq 0, 1\leq i\leq m$,

$$f(C_1,x_1)\leq f(C_2,x_2)\leq \cdots \leq f(C_m,x_m).$$

\no Otherwise, they will intersect.

$\mathcal{C}_{1100}$: In this case, a family ${\mathcal L}$ of
parallel lines is a parallel bundle along $\overrightarrow{O}$
only if there is an integer $k, 2\leq k\leq m$, such that

$$f(C_1)\leq f(C_2)\leq\cdots\leq f(C_{k-1})\ < \ f(C_k)=f(C_{k+1})=\cdots =f(C_m)=\pi.$$

\no Otherwise, they will intersect.

$\mathcal{C}_{1010}$: Similar to the case $\mathcal{C}_{1100}$, in
this case, a family ${\mathcal L}$ of parallel lines is a parallel
bundle along $\overrightarrow{O}$ only if there is an integer $k,
2\leq k\leq m$, such that

$$\pi =f(C_1)= f(C_2)=\cdots = f(C_k) \ < \ f(C_{k+1})\leq\cdots\leq f(C_m).$$

\no Otherwise, they will intersect.

$\mathcal{C}_{1001}$: In this case, a family ${\mathcal L}$ of
parallel lines is a parallel bundle along $\overrightarrow{O}$
only if there is an integer $k, l, 1\leq k \ < \ l\leq m$, such
that for $\forall x_i, x_i\geq 0, 1\leq i\leq k$ or $l\leq i \leq
m$,

\begin{eqnarray*}
f(C_1,x_1)&\leq& f(C_2,x_2)\leq\cdots\leq f(C_k,x_k)\ <
\ f(C_{k+1})\\
&=& f(C_{k+2})=\cdots =f(C_{l-1})=\pi \ < \
f(C_l,x_l)\leq\cdots\leq f(C_m,x_m).
\end{eqnarray*}

\no Otherwise, they will intersect.

$\mathcal{C}_{0110}$: In this case, a family ${\mathcal L}$ of
parallel lines is a parallel bundle along $\overrightarrow{O}$
only if there is integers $k, 1\leq k \ < \ m$, such that

$$f(C_1)\leq f(C_2)\leq\cdots\leq f(C_k) \ < \ \pi \ < \ f(C_{k+1})\leq\cdots\leq f(C_m).$$

\no Otherwise, they will intersect.

$\mathcal{C}_{0101}$: In this case, a family ${\mathcal L}$ of
parallel lines is a parallel bundle along $\overrightarrow{O}$
only if there is integers $k, 1\leq k \leq m$, such that for
$\forall x_i, x_i\geq 0, 1\leq i\leq m$,

$$ f(C_1,x_1)\leq f(C_2,x_2)\leq\cdots\leq f(C_k,x_k) \ < \pi \ \leq f(C_{k+1},x_{k+1})
\leq \cdots\leq f(C_m,x_m),$$

\no and there must be a constant in $f(C_1,x_1), f(C_2,x_2),\cdots
,f(C_k,x_k)$.

$\mathcal{C}_{0011}$: In this case, the situation is similar to
the case $\mathcal{C}_{0101}$ and there must be a constant in
$f(C_{k+1},x_{k+1}), f(C_{k+2},x_{k+2}),\cdots ,f(C_m,x_m)$.

$\mathcal{C}_{1110}$: In this case, a family ${\mathcal L}$ of
parallel lines is a parallel bundle along $\overrightarrow{O}$
only if there is an integer $k, l, 1\leq k \ < \ l\leq m$, such
that

\begin{eqnarray*}
f(C_1)&\leq& f(C_2)\leq\cdots\leq f(C_k) \ < \ f(C_{k+1})\\
&=&\cdots =f(C_{l-1})=\pi \ < \ f(C_l)\leq\cdots\leq f(C_m).
\end{eqnarray*}

\no Otherwise, they will intersect.

$\mathcal{C}_{1101}$: In this case, a family ${\mathcal L}$ of
parallel lines is a parallel bundle along $\overrightarrow{O}$
only if there is an integer $k, l, 1\leq k \ < \ l\leq m$, such
that for $\forall x_i, x_i\geq 0, 1\leq i\leq k$ or $l\leq i\leq
m$,

\begin{eqnarray*}
f(C_1,x_1)&\leq & f(C_2,x_2)\leq\cdots\leq f(C_k,x_k) \ < \
f(C_{k+1})\\
&=&\cdots =f(C_{l-1})=\pi \ < \ f(C_l,x_l)\leq\cdots\leq
f(C_m,x_m)
\end{eqnarray*}

\no and there must be a constant in $f(C_1,x_1), f(C_2,x_2),\cdots
,f(C_k,x_k)$. Otherwise, they will intersect.

$\mathcal{C}_{1011}$: In this case, a family ${\mathcal L}$ of
parallel lines is a parallel bundle along $\overrightarrow{O}$
only if there is an integer $k, l, 1\leq k \ < \ l\leq m$, such
that for $\forall x_i, x_i\geq 0, 1\leq i\leq k$ or $l\leq i\leq
m$,

\begin{eqnarray*}
f(C_1,x_1)&\leq & f(C_2,x_2)\leq\cdots\leq f(C_k,x_k) \ < \
f(C_{k+1})\\
&=&\cdots =f(C_{l-1})=\pi \ < \ f(C_l,x_l)\leq\cdots\leq
f(C_m,x_m)
\end{eqnarray*}

\no and there must be a constant in $f(C_l,x_l),
f(C_{l+1},x_{l+1}),\cdots ,f(C_m,x_m)$. Otherwise, they will
intersect.

$\mathcal{C}_{0111}$: In this case, a family ${\mathcal L}$ of
parallel lines is a parallel bundle along $\overrightarrow{O}$
only if there is an integer $k, 1\leq k \leq m$, such that for
$\forall x_i, x_i\geq 0$,

$$f(C_1,x_1)\leq f(C_2,x_2)\leq\cdots\leq f(C_k,x_k) \ < \ \pi
\ < \ f(C_l,x_l)\leq\cdots\leq f(C_m,x_m)$$

\no and there must be a constant in $f(C_1,x_1), f(C_2,x_2),\cdots
,f(C_k,x_k)$ and a constant in $f(C_l,x_l),
f(C_{l+1},x_{l+1}),\cdots ,f(C_m,x_m)$. Otherwise, they will
intersect.

$\mathcal{C}_{1111}$: In this case, a family ${\mathcal L}$ of
parallel lines is a parallel bundle along $\overrightarrow{O}$
only if there is an integer $k, l, 1\leq k \ < \ l\leq m$, such
that for $\forall x_i, x_i\geq 0, 1\leq i\leq k$ or $l\leq i\leq
m$,

\begin{eqnarray*}
f(C_1,x_1)&\leq & f(C_2,x_2)\leq\cdots\leq f(C_k,x_k) \ < \
f(C_{k+1})\\
&=&\cdots =f(C_{l-1})=\pi \ < \ f(C_l,x_l)\leq\cdots\leq
f(C_m,x_m)\end{eqnarray*}

\no and there must be a constant in $f(C_1,x_1), f(C_2,x_2),\cdots
,f(C_k,x_k)$ and a constant in $f(C_l,x_l),
f(C_{l+1},x_{l+1}),\cdots ,f(C_m,x_m)$. Otherwise, they will
intersect.

Following the structural characters of the classes
$\mathcal{C}_{1000}\sim\mathcal{C}_{1111}$, by the Proposition
$2.2$, $2.3$ and Proposition $3.1$, we know that any parallel
bundle is in one of the classes ${\mathcal C}_{1000}\sim{\mathcal
C}_{1111}$ and each class in ${\mathcal C}_{1000}\sim{\mathcal
C}_{1111}$ is non-empty. This completes the proof.\quad\quad
$\natural$

A example of parallel bundle in a planar map geometry is shown in
the Fig.$5$, in where the number on a vertex $u$ denotes the
number $\rho (u)\mu (u)$.

\includegraphics[bb=80 0 300 250]{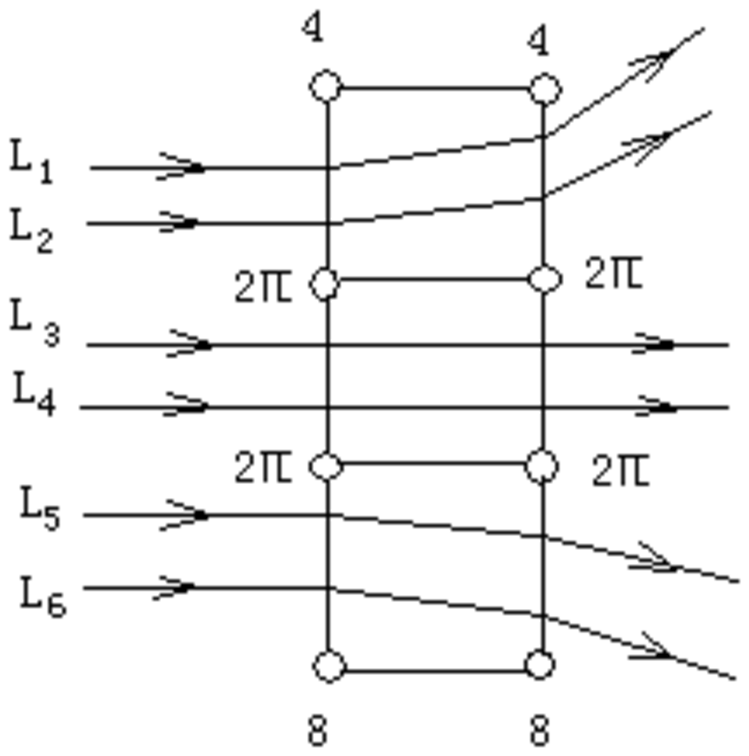}

\vskip 2mm

\c{\bf Fig.$5$}

\vskip 8mm

{\bf $5.$ Generalization}

\vskip 3mm

All the planar map geometries considered in this paper are without
boundary. For planar map geometries with boundary, i.e., some
faces are deleted ($[10]$), which are correspondence with the maps
with boundary ($[2]$). We know that they are the Smarandache {\it
non-geometries}, satisfying one or more of the following
conditions:

\vskip 3mm

($A1^-$){\it It is not always possible to draw a line from an
arbitrary point to another arbitrary point.}

($A2^-$){\it It is not always possible to extend by continuity a
finite line to an infinite line.}

($A3^-$){\it It is not always possible to draw a circle from an
arbitrary point and of an arbitrary interval.}

($A4^-$){\it not all the right angles are congruent.}

($A5^-$){\it if a line, cutting two other lines, forms the
interior angles of the same side of it strictly less than two
right angle, then not always the two lines extended towards
infinite cut each other in the side where the angles are strictly
less than two right angle.}\vskip 2mm

Notice that for an one face planar map geometry $(M,\mu)^{-1}$
with boundary, if we choose all points being euclidean, then
$(M,\mu)^{-1}$ is just the Poincar\'{e}'s model for the hyperbolic
geometry.

Using the neutrosophic logic idea, we can also define the
conception of {\it neutrosophic surface} as follow, comparing also
with the surfaces in $[8]$ and $[14]$.

\vskip 3mm

\no{\bf Definition $5.1$} {\it A neutrosophic surface is a
Hausdorff, connected, topological space $S$ such that every point
$v$ is elleptic, euclidean, or hyperbolic. }

\vskip 2mm

For this kind of surface, we present the following problem for the
further researching.

\vskip 3mm

\no{\bf Problem $5.1$} {\it To determine the behaviors of
elements, such as, the line, angle, area, $\cdots$, in
neutrosophic surfaces.}\vskip 2mm

Notice that results in this paper are just the behaviors of line
bundles in a neutrosophic plane.

\vskip 8mm

{\bf References}\vskip 3mm

\re{[1]}A.D.Aleksandrov and V.A.Zalgaller, {\it Intrinsic geometry
of surfaces}, American Mathematical Society, 1967.

\re{[2]}R.P.Bryant and D.Singerman, Foundations of the theory of
maps on surfaces with boundary,{\it
Quart.J.Math.Oxford}(2),36(1985), 17-41.

\re{[3]}H.Iseri, {\it Smarandache manifolds}, American Research
Press, Rehoboth, NM,2002.

\re{[4]}H.Iseri, {\it Partially Paradoxist Smarandache
Geometries}, http://www.gallup.unm.
edu/\~smarandache/Howard-Iseri-paper.htm.

\re{[5]}L.Kuciuk and M.Antholy, An Introduction to Smarandache
Geometries, {\it Mathematics Magazine, Aurora, Canada},
Vol.12(2003)

\re{[6]}Y.P.Liu, {\it Enumerative Theory of Maps}, Kluwer Academic
Publisher, Dordrecht / Boston / London (1999).

\re{[7]}Y.P.Liu, {\it Embeddability in Graphs}, Kluwer Academic
Publisher, Dordrecht / Boston / London (1995).

\re{[8]}Mantredo P.de Carmao, {\it Differential Geometry of Curves
and Surfaces}, Pearson Education asia Ltd (2004).

\re{[9]}L.F.Mao, {\it Automorphism groups of maps, surfaces and
Smarandache geometries}, American Research Press, Rehoboth,
NM,2005.

\re{[10]}L.F.Mao, A new view of combinatorial maps by
Smarandache's notion, arXiv: Math.GM/0506232, will also appear in
{\it Smarandache Notion J.}

\re{[11]}B.Mohar and C.Thomassen, {\it Graphs on Surfaces}, The
Johns Hopkins University Press, London, 2001.

\re{[12]}V.V.Nikulin and I.R.Shafarevlch, {\it Geometries and
Groups}, Springer-Verlag Berlin Heidelberg (1987)

\re{[13]}F. Smarandache, Mixed noneuclidean geometries, {\it
eprint arXiv: math/0010119}, 10/2000.

\re{[14]}J.Stillwell, {\it Classical topology and combinatorial
group theory}, Springer-Verlag New York Inc., (1980).

\end{document}